# Une Nouvelle Condition d'indépendance pour le théorème de la limite centrale.

**René Blacher**
**Laboratoire LMC**
**BP 53**
**38.041 Grenoble Cedex 9 FRANCE**

**Résumé:** On prouve un théorème de la limite centrale avec des conditions d'indépendance asymptotique beaucoup plus faibles que les hypothèses classiques.

**Summary :** We prove a Central Limit theorem with assumptions which are many much weak than the classical assumptions.

## I Introduction

**Notations 1-1 :** Let $\{X_n\}$ be a sequence of random variables defined on a probability space $(\Omega,\mathcal{A},P)$ such that, for all n $\mathbb{N}$, $\mathbb{E}\{X_n\} = 0$ and $0 < \mathbb{E}\{|X_n|^2\} < \infty$ .

Let $\sigma(n)^2$ be the variance of $X_1+X_2+\ldots+X_n$ . One sets $S_n = (X_1+X_2+\ldots+X_n)/\sigma(n)$ .

---

There is two classical ways to study the central limit theorem (CLT). One can use mixing conditions or martingales.

**Strong Mixing condition 1-2** : The process $\{X_n\}$ is said to be strongly mixing with coefficient $\alpha$ if

$$\sup_{A \in \mathcal{M}_1^n,\, B \in \mathcal{M}_{n+h}^\infty} |P(A \cap B)-P(A)P(B)| = \alpha(h) \to 0 \text{ as } h \to \infty \quad \text{(I-1)},$$

where for a<b, $\mathcal{M}^b_a$ is the σ-field generated by $(X_a,X_{a+1},\ldots,X_b)$ .

---

In 1962 Ibragimov has obtained a necessary and sufficient condition under which strongly mixing sequences satisfy the CLT (cf [1] and Denker [11] p 269-274) ).

**Theorem 1-1** : Assume that $\{X_n\}$ is stricly stationary such that (I-1) holds. Assume that $\sigma(n) \to \infty$ as $n \to \infty$ .

Then, $S_n \xrightarrow{d} N(0,1)$ if and only if $(S_n)^2$ is uniformly integrable.

In this case, $\sigma^2(n) = nh(n)$ where h is a slowly varying function

---

Ibragimov has deduced some results for functionals of mixing sequences. For example, Ibragimov has proved the theorem 18-6-2 of [13]..

These results have been developed by several authors , e.g. Bradley [2]-[3], Peligrad [4], Dehling-Denker-Philips [5], Davydov [6], Chan [7], Tjøsteim [8], Liebscher [9], Utev [10], and many others , Mervelede-Peligrad [45], Neumeyer [46], Dedecker-Rio [47], Johnson Barron [48] ((cf also Eberlein-Taqq [11], Hall-Heyde [12)], Ibragimov -Linnik [13] and Doukhan [14] and [49]-[54]).

The CLT has been proved for martingales also cf [12]. Theses results have been developed by Chao [40], De Meyer [41], Wang-Yang-Zhou [42], Mervelede [43], Ouchti [44]).

In the most part of these papers, the CLT is studied with the convergence in distribution. However, some authors have studied the moment’s convergence (Bernstein [19], Brown [20], Eissein-Janson [21], Hernndoff [22], Birkel [23], Krugov [24] , Mairoboda [25], Yokohama [26],[27] , Ibragimov [55], Soulier [56], Rozovsky [57]). For example, Yokohama has obtained necessary conditions in order to the moments converges. Cox-Kim have studied the moments bounds [28].

Then, under some assumptions (cf [27] and [34], and [12] p 71) the MCLT holds (Moments’s Central Limit theorem) : $S_n \xrightarrow{M} N(0,1)$ , that is, for all $p \in \mathbb{N}$, $\mathbb{E}\{(S_n)^p\} \to \mu_p$ as $n \to \infty$ , where $\mu_p$ is the p-th moment of $N(0,1)$ .

Now, mixing condition or martingale condition are the most used assumptions. But, theses conditions are not necessary conditions. As a mater of fact, those are strong assumptions : e.g. Strong mixing condition does not hold for some AR(1) processes (cf [13] p 360-362, [11] p 180). However the CLT holds (cf th 18-6-5 of [13])

Then, classical conditions are too strong. As a matter of fact, those are a tiny minority of asymptotic independence conditions which are sufficient on order that the CLT holds .

Because strong mixing condition is too strong and martingale condition is too specific, some authors have introduce weaker hypotheses : Versik Ornstein (cf [2] [5]) , Withers [15], Cogburn [16] (cf also Rosenblatt [17]), Pinsker [18] Doukhan-Prieur [62].

Theses conditions are more general than mixing conditions or martingale conditions but they are not necessary conditions. As matter of fact they are not founded on a measure of dependence which determines completely dependence. Then, in order to obtain a full solution, we have used a such measure : the higher order correlation coefficients $\rho_{j_1,j_2,....,j_n}$. Theses coefficients are perfectly adapted to MCLT.

Because dependence is completely determined by the $\rho_{j_1,j_2,....,j_n}$ 's , one can ask what is their part in the MCLT. In [33] one has obtained necessary and sufficient conditions for MCLT . One has deduced almost minimal assumptions for the MCLT.

## II Generalization of mixing conditions for MCLT

Now, this type of assumptions is different of 1.2. Then, one wants a condition whose the writing is near of that of mixing condition.

With this aim, one decomposes $X_1+X_2+.....+X_n$ in $X_1+X_2+.....+X_u$ , $X_{u+1}+X_{u+2}+.....+X_{u+t}$ and $X_{u+t+1}+X_{u+t+2}+.....+X_{u+t+u}$ , where u and t are functions of n : u = u(n), t=t(n), u+t+u = n.

**Notations 2-1 :** We denote by $\kappa(n) \in \mathbb{N}$, an increasing sequence such that $\kappa(1)= 0$, $\kappa(n)\leq n$ and $\kappa(n)/n \to 0$ as $n \to \infty$ . We define the sequences u=u(n) and t=t(n) by : u(1)=1, $u(n) = \max\{m \in \mathbb{N}^* \mid 2m+\kappa(m) \leq n\}$ and t(1)=0, t(n) = n-2u(n) if $n\geq 2$.

---

We shall choose $\kappa$ such that $u(n)\to\infty$, $n/u(n) \to 2$ and $t(n)/u(n) \to 0$ as $n\to\infty$ (cf 5-18). Moreover, if $\kappa(n)\to\infty$, $t(n)\to\infty$ as $n\to\infty$.

Now, we need to normalize the associated sums. Then, one tolerates the following notations (cf 1-1).

**Notations 2-2 :** Let $\sigma(u)^2$ be the variance of $X_1+X_2+...+X_u$ . One sets $S_u = (X_1+X_2+...+X_u)/\sigma(u)$, $\xi_u = (X_{u+1}+X_{u+2}+.....+X_{u+t})/\sigma(u)$ and $S'_u = (X_{u+t+1}+X_{u+t+2}+.....+X_{u+t+u})/\sigma(u)$. If t(n)=0, we set $\xi_u = 0$.

---

Then, one can assume $\mathbb{E}\{(S_u)^p(S'_u)^q\} - \mathbb{E}\{(S_u)^p\}\mathbb{E}\{(S'_u)^q\}\to 0$ for all (p,q) as
M
asymptotic independence condition. Indeed, if $\xi_u$ is negligible $S_n \to \ N(0,1)$ .

Now in order to prove the MCLT for some distributions, one generalize this condition by the following way (cf [61]).

**Theorem 2-1** : Assume $\mathbb{E}\{X_n\}= 0$ for all n $\mathbb{N}^*$. Let k $\mathbb{N}^*$. Assume that, for all p $\mathbb{N}$, p<2k+1, $\mathbb{E}\{(X_n)^p\}< \infty$ for all n $\mathbb{N}^*$ and that $\mathbb{E}\{(\xi_u)^p\}\to 0$ as $n \to \infty$ .

Assume that there exists two sequences of random variables $\{\upsilon'_u\}$ and $\{\upsilon'_u\}$, $\mathbb{E}\{(\upsilon_u)^p\}+\mathbb{E}\{(\upsilon'_u)^p\} \to 0$ as $n \to \infty$ for all p $\mathbb{N}$, p<2k+1, such that the following assumptions hold.

**$H_{mS}$(2k)** : $\forall$p $\mathbb{N}$, p<2k+1, $\mathbb{E}\{(S_u+\upsilon_u)^p\} - \mathbb{E}\{(S'_u+\upsilon'_u)^p\}\to 0$ as $n\to\infty$ .

**$H_{mI}$(2k) :** $\forall$(p,q) $\mathbb{N}^{*2}$, p+q<2k+1,

$$\mathbb{E}\{(S_u+\upsilon_u)^p(S'_u+\upsilon'_u)^q\} - \mathbb{E}\{(S_u+\upsilon_u)^p\}\mathbb{E}\{(S'_u+\upsilon'_u)^q\}\to 0 \text{ as } n\to\infty .$$

Then, for all p $\mathbb{N}$, p≤2k, $\mathbb{E}\{(S_n)^p\} \to \mu_p$ .

---

If k=∞, we set $H_{mS}(\infty) = H_{mS}$ and $H_{mI}(\infty)=H_{mI}$ . Then, compare $H_{mI}$ and $H_{mS}$ with classical assumptions. Of course, if $\{X_n\}$ is strictly stationary, $H_{mS}$ holds.

Assume that $\{X_n\}$ is a strictly stationary φ-mixing process. By theorem 17-2-3 of [13],

$$\left|\mathbb{E}\{(S_u)^p(S'_u)^q\} - \mathbb{E}\{(S_u)^p\}\mathbb{E}\{(S'_u)^q\}\right| \le 2\varphi(t)\ \mathbb{E}\{|S_u|^{pa}\}^{1/a}\ \mathbb{E}\{|S_u|^{qb}\}^{1/b} ,$$

where a,b>1 $a^{-1}+b^{-1} = 1$.

Similarly, if $\{X_n\}$ is a strictly stationary strongly mixing process, by theorem 17-2-2 of [13], there exists a function $\mathfrak{F}$ : $\mathbb{R}^3\to \mathbb{R}_+$ such that, for all (p,q) $\mathbb{N}^{*2}$, for all n $\mathbb{N}^*$,

$$\mathfrak{F}(\beta,p,q)\left|\mathbb{E}\{(S_u)^p(S'_u)^q\} - \mathbb{E}\{(S_u)^p\}\mathbb{E}\{(S'_u)^q\}\right| \le \alpha(t)^{1-\beta} ,$$

where β≥0.

Remark that it is not needed that $H_{mI}$ holds if $\{X_n\}$ is strong mixing. But it is the case
M
if $S_n \to \ N(0,1)$. In this case $H_{mI}$ is weaker than the strong mixing assumption : no conditions are dictated to the rate of convergence of $\mathbb{E}\{(S_u)^p(S'_u)^q\} - \mathbb{E}\{(S_u)^p\}\mathbb{E}\{(S'_u)^q\}$.

**Example 2-3 :** For example, if $\Theta_1$ has the standart normal distribution and $g_i(\Theta_1) = i^{-3}\sin(2\pi i\Theta_1)$, $S_n \xrightarrow{M} N(0,1)$ (cf [61]).

**Example 2-4 :** Let $\{\Theta_t\}$ be an IID sequence of random variables such that $\Theta_1$ has the uniform distribution on [-1,1]. Let $\{\Psi_t\}$ be a strictly stationary process independent of $\{\Theta_t\}$. Let $X_t$ be the process

$$X_t = \sum_{i\geq 0} i^{-3/2} L_{2^i}(\Theta_{t+i}) f_i(\Psi_{t+i}) \quad ,$$

where $\{L_i\}$ is the family of Legendre polynomials such that $L_i(x) = x^i + b_{i-1}x^{i-1} + \ldots\ldots + b_0$ , and where $|f_i(\Psi_1)| \leq 1$. One assumes that $\{X_t\}$ satisfies the first assumptions of theorem 1-7 of [33].

Then, $\mathbb{E}\{L_i(\Theta_1)\pi(\Theta_1)\} = 0$ for all polynomial $\pi$ such that $\deg(\pi) < i$. One deduces that $X_t$ satisfies the second assumptions of theorem 1-7 of [33]., but is not a martingale.

Then, by theorem 1-7 of [33], $S_n \xrightarrow{M} N(0,1)$ .

One deduce that $H_{mI}$ holds.

Remark that the dependence between $\{X_t\}$ and $\{X_{t+h}\}$ can be strong, e.g., if $f_i(\Psi_t) = \Psi_t$, it is enough that $\left|\mathbb{E}\{(\Psi_t)^2(\Psi_{t+h})^2\} - \mathbb{E}\{(\Psi_t)^2\}\mathbb{E}\{(\Psi_{t+h})^2\}\right| \leq e(h) \rightarrow 0$ in order that $S_n \xrightarrow{M} N(0,1)$ .

# III Application to CLT

Theorem 2-1 suggests that one can replace the strong mixing condition by the following way.

**Definition 3-1** : We define condition $\mathbf{H_S}$ and $\mathbf{H_I}$ by the following way. There exists two sequences of random variables $\upsilon_u$ and $\upsilon'_u$, $\mathbb{E}\{\upsilon_u\} = \mathbb{E}\{\upsilon'_u\} = 0$, $\mathbb{E}\{(\upsilon_u)^2\} + \mathbb{E}\{(\upsilon'_u)^2\} \rightarrow 0$ as $n \rightarrow \infty$ , such that

$\mathbf{H_S}$ : $\forall k \ \mathbb{N}$, $\forall j \ \mathbb{N}$, $\mathbf{P}\{A_{k,j}\} - \mathbf{P}\{B_{k,j}\} \rightarrow 0$ as $n \rightarrow \infty$ ,

$\mathbf{H_I}$: $\forall k \ \mathbb{N}$, $\forall (j,j') \ \mathbb{N}^2$ , $\mathbf{P}\{A_{k,j} \cap B_{k,j'}\} - \mathbf{P}\{A_{k,j}\}\mathbf{P}\{B_{k,j'}\} \rightarrow 0$ as $n \rightarrow \infty$ ,

where $A_{k,j}$ and $B_{k,j}$ are the events $A_{k,j} = \{\sigma^*(u)^{-1}(S_u+\upsilon_u)\ \ [4^{-k}j,4^{-k}(j+1)[\}$ and $B_{k,j} = \{\sigma^*(u)^{-1}(S'_u+\upsilon'_u)\ \ [4^{-k}j,4^{-k}(j+1)[\}$ with $\sigma^*(u)^2 = \mathbb{E}\{(S_u+\upsilon_u)^2\}$.

---

In particular, if $\upsilon_u = \upsilon'_u = 0$, $A_{k,j} = \{S_u\ \ [4^{-k}j\ ,\ 4^{-k}(j+1)[\ \}$ and $B_{k,j} = \{S'_u\ \ [4^{-k}j\ ,\ 4^{-k}(j+1)[\}$ . Remark also that $\sigma^*(u)^2 \to 1$.

Moreover, if the CLT holds and if $\mathbf{H_S}$ and $\mathbf{H_I}$ hold for $\upsilon_u = \upsilon'_u = 0$, then, $\mathbf{H_S}$ and $\mathbf{H_I}$ hold also for all sequences $\upsilon_u$ and $\upsilon'_u$ such that $\mathbb{E}\{\upsilon_u^2\}+\mathbb{E}\{\upsilon'^2_u\}\to 0$.

Now, generally, it is complicated to prove that $H_I$ holds for $\upsilon_u = \upsilon'_u = 0$. Then it is simpler to use the above assumptions. For example, if $\Psi_t = \theta_t$ , $\mathbf{H_I}$ and $\mathbf{H_S}$ hold if $t(n)\to\infty$ with $\upsilon_u = -r_u$ and $\upsilon'_u = -r'_u$ .

Then the following theorem holds.

**Theorem 3-1** : We keep the previous notations. Assume that $\mathbf{H_S}$ and $\mathbf{H_I}$ hold. Assume that $\mathbb{E}\{S_u^2\} - \mathbb{E}\{S'^2_u\}\to 0$ and $\mathbb{E}\{\xi_u^2\}\to 0$ as $n\to\infty$.

Then, $S_n \xrightarrow{d} N(0,1)$ if and only if $(S_n+\upsilon_n)^2$ is uniformly integrable .

In this case, $\sigma(n)\to\infty$ and $\sigma(n)/\sigma(u)\to 2^{1/2}$ as $n\to\infty$ .

---

Remark that $n/u \to 2$ (cf lemma 5-18). Then, "$\sigma(n)/\sigma(u)\to 2^{1/2}$" is weaker than "h is a slowly varying function", that is $\sigma(tn)^2/\sigma(n)^2\to t$" for all t>0 (cf [13], 18-2-3, p 325 and 394).

One can also compare $\mathbf{H_I}$ with the other generalizations of the strong mixing condition : [2],[5],[15],[16], [17] , [18] and [62] .

Of course, in $\mathbf{H_I}$ the uniformity of the strong mixing condition is suppressed. Moreover, if $\{X_n\}$ is strongly mixing, $\mathbf{H_I}$ holds : $\sup\left|\mathbf{P}\{A_{k,j}\cap B_{k,j}\}-\mathbf{P}\{A_{k,j}\}\mathbf{P}\{B_{k,j}\}\right| \le \alpha\{t(n))$.

Now, in some cases, it may be simpler to use $\mathbf{H_{mI}}$ than $\mathbf{H_I}$ .

**Corollary 3-2** : Assume that $\mathbb{E}\{|X_n|^4\}<\infty$ and that $\mathbb{E}\{(\xi_u)^4\}\to 0$ . Assume that $\mathbf{H_S}$ , $\mathbf{H_I}$ , $\mathbf{H_{mS}(4)}$ and $\mathbf{H_{mI}(4)}$ hold.

Then , $S_n \xrightarrow{d} \mathcal{N}(0,1)$. Moreover, $\mathbb{E}\{(S_n)^3\} \to 0$ and $\mathbb{E}\{(S_n)^4\} \to 3$.

---

**Proposition 3-3** : Assume that the assumptions of theorem 2-1 hold with k=∞. Then, $\mathbf{H_I}$ and $\mathbf{H_S}$ hold.

---

**Proof** : By our assumptions, $\mathbb{E}\{(S_u+\upsilon_u)^k\} \to \mu_k$ . By $\mathbf{H_{mS}}$ , $\mathbb{E}\{(S'_u+\upsilon'_u)^k\} \to \mu_k$ . By $\mathbf{H_{mI}}$ , for all (a,b) $\mathbb{R}^2$, $\mathbb{E}\Big\{\big(a(S_u+\upsilon_u)+b(S'_u+\upsilon'_u)\big)^k\Big\}$ converges to the k-th moment of $\mathcal{N}(0,a^2+b^2)$. One deduces that $(S_u+\upsilon_u\, ,S'_u+\upsilon'_u) \xrightarrow{d} \mathcal{N}_2(0,\mathbf{I}_2) = \mathcal{N}(0,1)\ \ \mathcal{N}(0,1)$.

---

## IV The central limit theorem for functionals.

Let $\{X_t\}$ be the process $X_t = \sum_{i\geq 0} g_{i+1}(\Theta_{t+i}) f_{i+1}(\Psi_{t+i})$ , where $\{\Psi_t\}$ is a strictly stationary process and where $\{\Theta_t\}$ is an IID sequence independent of $\{\Psi_t\}$ . We assume that, for all i $\mathbb{N}$, $g_i$ and $f_i$ are measurable, $\mathbb{E}\{g_i(\Theta_1)f_i(\Psi_1)\} = 0$ and $\mathbb{E}\Big\{g_i(\Theta_1)^2 f_i(\Psi_1)^2\Big\} <\infty$ .

Then, one can write $S_n=M_n+r_n$ with $M_n=\frac{1}{\sigma(n)}\sum_{r=1}^{n}\sum_{s=1}^{r} g_s(\Theta_r)f_s(\Psi_r)$ , $r_n= \frac{1}{\sigma(n)}\sum_{r=1}^{\infty}\sum_{s=r+1}^{r+n} g_s(\Theta_{n+r})f_s(\Psi_{n+r})$ : $\mathbb{E}\{M_n\}= \mathbb{E}\{r_n\}=0$ .

Of course, $S_{u(n)} = M_{u(n)} + r_{u(n)}$ . Moreover, one can write $S'_{u(n)} = M'_{u(n)} + r'_{u(n)}$ with $M'_u = \frac{1}{\sigma(u)}\sum_{r=1}^{n}\sum_{s=1}^{r} g_s(\Theta_{r+u+t})f_s(\Psi_{r+u+t})$ and $r'_u = \frac{1}{\sigma(u)}\sum_{r=1}^{\infty}\sum_{s=r+1}^{r+n} g_s(\Theta_{n+r+u+t})f_s(\Psi_{n+r+u+t})$.

Under some assumptions, one can proves that $H_I$ holds. In particular that is the case if $\Psi_t = \theta_t$ where $\{\theta_t\}$ is a strongly mixing process. Indeed, it is easy to prove the following proposition.

**Proposition 4-1** : Assume that $X_t = \sum_{i\geq 0} g_{i+1}(\Theta_{t+i}) f_{i+1}(\theta_{t+i})$ with $\mathbb{E}\{g_i(\Theta_1)\} = 0$ for all $i \in \mathbb{N}$. Assume that $\sum_{1\leq i\leq r} g_i(\Theta_1) f_i(\theta_1)$ converges in $L^2(\Omega)$ to $G$ such that $\int G^2.dP > 0$ with $\sum_{s=1}^{n} \int (\sum_{i=s}^{\infty} g_i(\Theta_1) f_i(\theta_1))^2.dP < C < \infty$ .

Then, $\sigma(n) \to \infty$ and $\mathbb{E}\{(r_n)^2\} \to 0$ as $n \to \infty$ . Moreover, $\mathbf{H_I}$ and $\mathbf{H_S}$ hold with $\upsilon_u = - r_u$ and $\upsilon'_u = - r'_u$ if $t(n) \to \infty$ .

---

Therefore, under the assumptions of 4-1, $\{X_n\}$ satisfies the CLT if and only if $(M_n)^2$ is uniformly integrable.

Now,we know that $(S_n + \upsilon_n)^2$ is uniformly integrable if $\mathbb{E}\{|S_n + \upsilon_n|^{2+\delta}\} \leq C < \infty$ where $\delta > 0$ ([11], p 270). Then, one can use the following property.

**Proposition 4-2** : Assume that the assumptions of 4-1 hold. We set $G_r(\Theta_1, \theta_1) = \sum_{s=1}^{r} g_s(\Theta_1) f_s(\theta_1)$ . Then, there exists K>0 such that $\mathbb{E}\{(M_n)^4\} < K$ if and only there exists B>0 such that $n^{-2} \sum_{r=1}^{n} \mathbb{E}\{G_r(\Theta_r, \theta_r)^4\} \leq B < \infty$.

---

Assume that $\mathrm{Sup} \left| \sum_{i=s}^{\infty} g_i(\Theta_1) f_i(\theta_1) \right| = \varepsilon(s)$ with $\sum_{s=1}^{\infty} \varepsilon(s) < C < \infty$ . Then, there exists $M^* > 0$ such that $|G_r(\Theta_r, \theta_r)| \leq M^*$ . Therefore, there exists K>0 such that $\mathbb{E}\{(M_n)^4\} \leq K$.

**4-3 Example** : Assume $g_i(\Theta_1) f_i(\theta_1) = (1/i)^3 \exp\{-i(\Theta_1)^2\} \cos(2\pi i \theta_1)$. Then, $S_n \xrightarrow{d} N(0,1)$ .

---

Then, in these two examples, $S_n \xrightarrow{d} N(0,1)$ . Remark that no condition on $\alpha(n)$, the mixing coefficient of $\{\theta_t\}$ is necessary.

**4-3 Generalization :** For example, let $X_t = \sum_{i\geq 0} f_{i+1}(\theta_{t+i})$ where $F_r(\theta_1) = \sum_{1\leq i\leq r} f_i(\theta_1)$ converges in $L^2(\Omega)$ to F such that $\int F^2.dP > 0$ . Let $e(s) = \int \left(\sum_{i=s}^{\infty} f_i(\theta_1)\right)^2.dP$ .

*If $\{\theta_t\}$ is strong mixing, one can assume $\sum_{s=1}^{\infty} e(s)^{1/2} < \infty$ in order that $H_I$ holds.*

# V : Generalization of Theorem 3-1.

## V-1 : Notations

In this section , one uses assumptions more general. Remark that we do not assume $\mathbb{E}\{(X_n)\} = 0$ and $\mathbb{E}\{(X_n)^2\} < \infty$ .

**Notations 5-1** : Let N(n) be an increasing sequence. Let $\{X_{n,t}\}$, n,t $\in \mathbb{N}^*$, $1\leq t\leq N(n)$, be a triangular array of random variables defined on $(\Omega, \mathcal{A}, \mathbf{P})$. For any r $\in \mathbb{N}^*$, we denote by $\Psi_r(n)$ a real sequence such that $\Psi_r(n)>0$.

For all n $\in \mathbb{N}^*$ and all r $\in \mathbb{N}^*$, such that $u^r(n)>1$, one sets

$$S_{r,n} = \Psi_r(n)^{-1}\left(X_{n,1}+X_{n,2}+...+X_{n,u^r(N)}\right) ,$$

$$\xi_{r,n} = \Psi_r(n)^{-1}\left(X_{n,u^r(N)+1}+X_{n,u^r(N)+2}+......+X_{n,u^r(N)+t_r(N)}\right) ,$$

$$S'_{r,n} = \Psi_r(n)^{-1}\left(X_{n,u^r(N)+t_r(N)+1}+X_{n,u^r(N)+t_r(N)+2}+......+X_{n,u^r(N)+t_r(N)+u^r(N)}\right) ,$$

where N= N(n), $t_r(n) = t\left(u^{r-1}(N)\right)$ , and $u^r(N)= u[u[....[u(N)]...]]$ .

We also write

$$S_n = S_{0,n} = \Psi_0(n)^{-1}\left(X_{n,1}+X_{n,1}+...+X_{n,N}\right) .$$

Let $\mathcal{J}_{k,j}$, $k \in \mathbb{N}$, be a sequence of nested partitions of $\mathbb{R}$ : $\mathbb{R} = \cup_{j \in \mathcal{N}_k} \mathcal{J}_{kj}$ where the $\mathcal{J}_{kj}$ 's are intervals and $\mathcal{N}_k \subset \mathbb{N}$. One assumes that, for all $k \in \mathbb{N}$, and all bounded interval $\mathcal{K}$, there exists a finite set $\mathcal{Z}_k \subset \mathcal{N}_k$ such that $\mathcal{K} \subset \cup_{j \in \mathcal{Z}_k} \mathcal{J}_{k,j}$ . One assumes also that $|\mathcal{J}_{kj}| \le e_k$ where $e_k \to 0$ as $k \to \infty$, where $|\mathcal{J}_{kj}|$ is the length of $\mathcal{J}_{kj}$.

Moreover, we write the uniform integrability of $(S_n + \upsilon_n)^2$ in the form :

$$\mathbb{E}_{|S_n + \upsilon_n| \ge k}\{(S_n + \upsilon_n)^2\} \le \Delta_k \text{ , with } \Delta_k \to 0 \text{ as } k \to \infty,$$

where $\mathbb{E}_{|S_n + \upsilon_n| \ge k}\{(S_n + \upsilon_n)^2\} = \int \P_{[-\infty,-k[\cup[k,\infty[}(S_n + \upsilon_n).\ (S_n + \upsilon_n)^2.dP$ .

---

As a matter of fact, we shall decompose $S_{0,n}$ in the form

$$(\Psi_0(n)/\ \Psi_1(n))\ S_{0,n} = \Psi_1(n)^{-1}\big(X_{n,1} + X_{n,1} + \ldots + X_{n,N}\big)$$

$$= S_{1,n} + \xi_{1,n} + S'_{1,n}$$

$$= \Psi_1(n)^{-1}\big[(X_{n,1} + \ldots + X_{n,u(N)}) + (X_{n,u(N)+1} + \ldots + X_{n,u(N)+t(N)}) + (X_{n,u(N)+t(N)+1} + \ldots + X_{n,n})\big] .$$

After, we shall decompose $S_{1,n}$ in

$$(\Psi_1(n)/\ \Psi_2(n))S_{1,n} = S_{2,n} + S'_{2,n} + \xi_{2,n}$$

$$= \Psi_2(n)^{-1}\big[(X_{n,1} + \ldots\ldots + X_{n,u(u(N))}) + (X_{n,u(u(N))+1} + \ldots.; + X_{n,u(u(N))+t(u(N))}) + (X_{n,u(u(N)))+t(u())+1} + \ldots.. + X_{n,n})\big].$$

Moreover we replace partitions $[4^{-k}j, 4^{-k}(j+1)[$ by more general partition $\mathcal{J}_{kj}$ which have the same useful properties. For example, if $\mathcal{J}_{kj} = [4^{-k}j, 4^{-k}(j+1)[$ , $|\mathcal{J}_{kj}| = e_k = 4^{-k} \to 0$ as $k \to \infty$. Moreover, for all bounded interval $\mathcal{K}$ there exists $P,Q \in \mathbb{Z}$ such that $\mathcal{K} \subset [P,Q[$ . Then, $[P,Q[ = \cup[4^{-k}j, 4^{-k}(j+1)[$ where $P4^k \le j < Q4^k$.

Then, we shall prove that theorem 3-1 holds for events $A_{k,j} = \{\sigma^*(u)^{-1}(S_u + \upsilon_u) \in \mathcal{J}_{kj}\}$ and $B_{k,j} = \{\sigma^*(u)^{-1}(S'_u + \upsilon'_u) \in \mathcal{J}_{kj}\}$ .

**Assumptions 5-2 :** We assume that $e_k = 4^{-k}$ . Let $\mathcal{K}_k$ be an increasing sequence of bounded intervals such that $\mathcal{K}_k = \cup_{j \in \mathcal{Z}_k} \mathcal{J}_{k,i}$ where $\mathcal{Z}_k \subset \mathcal{N}_k$ is a finite set for all $k \in \mathbb{N}$ and $\mathbb{R} = \cup_{k \in \mathbb{N}} \mathcal{K}_k$ . We set $b_k = \sup\{|x| \mid x \in \mathcal{K}_k\}$ . We denote by $\mathcal{C}_k$ the $\sigma$-algebra generated by $\{\mathcal{J}_{k,j}\}_{j \in \mathcal{Z}_k}$ .

Let $k,h \in \mathbb{N}$ and $\tau>0$. One sets

$$\varepsilon^{I}_{k,h}(n) = \sup_{n'\geq n} \left\{ \sup_{\mathcal{B},\mathcal{B}' \in \mathcal{C}_{k},\ r\leq h} \left\{ \left| \mathbf{P}\{(S_{r,n'} \in \mathcal{B}) \cap (S'_{r,n'} \in \mathcal{B}')\} - \mathbf{P}\{S_{r,n'} \in \mathcal{B})\}\mathbf{P}\{S'_{r,n'} \in \mathcal{B}'\} \right| \right\} \right\} ,$$

$$\varepsilon^{S}_{k,h}(n) = \sup_{n'\geq n} \left\{ \sup_{\mathcal{B} \in \mathcal{C}_{k},\ r\leq h} \left\{ \left| \mathbf{P}\{S_{r,n'} \in \mathcal{B}\} - \mathbf{P}\{S'_{r,n'} \in \mathcal{B}'\} \right| \right\} \right\} ,$$

$$\varepsilon^{\psi}_{\tau,h}(n) = \sup_{n'\geq n,\ 0<r\leq h} \left\{ 1 - \frac{\psi_{r-1}\{N(n')\}}{2^{\tau}\psi_{r}\{N(n')\}} \right\} .$$

Moreover, one denotes by $\varepsilon^{\xi}_{h}(n)$ a non increasing sequence such that $\mathbf{P}\{|\xi_{r,n}|>\varepsilon^{\xi}_{h}(n)\} < \varepsilon^{\xi}_{h}(n)$ for all $r=1,2,...,h$.

One sets $\varepsilon^{\tau}_{k}(n) = b_{2k}\varepsilon^{\psi}_{\tau,k}(n)+\varepsilon^{\xi}_{k}(n) +\varepsilon^{I}_{2k,k}(n) +\varepsilon^{S}_{2k,k}(n)$ .

---

In paragraph V-3, we shall suppose that $\mathbf{P}\{(S_{r,n} \in \mathcal{J}_{k,j}) \cap (S'_{r,n} \in \mathcal{J}_{k,j'})\} - \mathbf{P}\{S_{r,n} \in \mathcal{J}_{k,j})\}\mathbf{P}\{S'_{r,n} \in \mathcal{J}_{k,j'}\} \to 0$ as $n\to\infty$ for all $k,j,j',r$. Then, because $\mathcal{C}_k$ is finite, $\varepsilon^{I}_{k,h}(n) \to 0$ as $n\to\infty$.

Moreover, we shall suppose also that, for all $r \in \mathbb{N}$, there exists a decreasing sequence $\varepsilon^{\xi}_{r}(n)$ such that $\mathbf{P}\{|\xi_{r,n}|>\varepsilon^{\xi}_{r}(n)\} < \varepsilon^{\xi}_{r}(n)$ for all $r=1,2,..,h$. Then, $\varepsilon^{\xi}_{h}(u)$ is also decreasing.

These assumptions more general as those of §III will can be used for a more complete study of the part of dependence coefficients in limit distributions (in particular, for the laws of large numbers and for the convergence to the Poisson distribution). We shall study these generalizations in full detail later.

Recall also that some CLT have been obtained for triangular array of random variables or if $\sigma(n)^2 \geq n^r$, r>2 (cf [8], [9], [10]). Non stationarity is studied in [7].

## V-2 : Lemmas

In this paragraph, one proves some inequalities under the previous hypotheses. Moreover, we assume that the following assumption holds.

**Hypothesis 5-3 :** In this paragraph V-2, we assume that, for all r=1,2,...,h,

$$\mathbf{P}\{S_{r,n} \;\mathfrak{K}_k)\} \leq 4^{-k} \quad \text{(V-1)}.$$

---

Then the following result holds.

**Lemma 5-4 :** For all r=1,2,...,h,

$$\mathbf{P}\{S'_{r,n} \;\mathfrak{K}_k)\} \leq 4^{-k} + \varepsilon^S_{k,h}(n) \quad \text{(V-2)}.$$

---

One proves the following inequalities by using the same way as Volkonskii and Rosanov ([38] : condition I'. cf also [37]).

**Lemma 5-5 :** Let k and h $\mathbb{N}$ and D>0. Let $\zeta = \sum_{j \;\mathfrak{J}_k} \alpha_j \P_{\mathfrak{J}_{k,j}}$ and $\eta = \sum_{j \;\mathfrak{J}_k} \beta_j \P_{\mathfrak{J}_{k,j}}$ where $|\alpha_j| \leq D$ and $|\beta_j| \leq D$. Then, for all r=1,2,....,h, the following inequalities hold :

$$\left|\mathbb{E}\{\zeta(S_{r,n})\} - \mathbb{E}\{\zeta(S'_{r,n})\}\right| \leq 2D\varepsilon^S_{k,h}(n)) \quad \text{(V-3)} \;,$$

$$\left|\mathbb{E}\{\zeta(S_{r,n})\eta(S'_{r,n})\} - \mathbb{E}\{\zeta(S_{r,n})\}E\{\eta(S'_{r,n})\}\right| \leq 4D^2\varepsilon^I_{k,h}(n) \quad \text{(V-4)} \;.$$

---

**Lemma 5-6 :** Let t $\mathbb{R}$ and k and h $\mathbb{N}$ . Then, for all r=1,2,...,h,

$$\left|\mathbb{E}\{e^{itS_{r,n}}\} - \mathbb{E}\{e^{itS'_{r,n}}\}\right| \leq 6\varepsilon^S_{k,h}(n) + 4\,.4^{-k} + 4|t|4^{-k} \quad \text{(V-5)} ,$$

$$\left|\mathbb{E}\{e^{itS_{r,n}+itS'_{r,n}}\} - \mathbb{E}\{e^{itS_{r,n}}\}\mathbb{E}\{e^{itS'_{r,n}}\}\right| \leq 16\varepsilon^I_{k,h}(n) + 8\varepsilon^S_{k,h}(n) + 16\,.4^{-k} + 16|t|4^{-k} \quad \text{(V-6)}.$$

---

**Proof :** There exists $C_{k;t} = \sum_{j \;\mathfrak{J}_k} \alpha_j^{k,t} \P_{\mathfrak{J}_{k,j}}$ , $|C_{k;t}| \leq 1$, such that $|\cos(ts) - C_{k;t}| \leq 4^{-k}|t|$ if s $\mathfrak{K}_k$'. Then, by (V-1), (V-2) and (V-3),

$$\left|\mathbb{E}\{\cos(tS_{r,n})\} - \mathbb{E}\{\cos(tS'_{r,n})\}\right|$$
$$\leq \left|\mathbb{E}_{S_{r,n}}\ \mathfrak{K}_k\{\cos(tS_{r,n})\} - \mathbb{E}_{S'_{r,n}}\ \mathfrak{K}_k\{\cos(tS'_{r,n})\}\right| + 2.4^{-k} + \mathcal{E}^{S}_{k,h}(n)$$
$$\leq \left|\mathbb{E}_{S_{r,n}}\ \mathfrak{K}_k\{C_{k;t}(S_{r,n})\} - \mathbb{E}_{S'_{r,n}}\ \mathfrak{K}_k\{C_{k;t}(S'_{r,n})\}\right| + 2.4^{-k} + 2.4^{-k}|t| + \mathcal{E}^{S}_{k,h}(n)$$
$$\leq 2.4^{-k} + 2.4^{-k}|t| + 3\mathcal{E}^{S}_{k,h}(n) \quad .$$

One uses the same way with $\mathbb{E}\{\sin(tS_{r,n})\}$, $\mathbb{E}\{\cos(tS_{r,n})\sin(tS_{r,n})\}$, etc (cf [37]).

---

**Lemma 5-7 :** Let t $\mathbb{R}$ and k and h $\mathbb{N}$ . Then, for all r=1,2,...,h,

$$\left|\mathbb{E}\{e^{it\psi_{r-1}(N).S_{r-1,n}/\psi_r(N)}\} - \mathbb{E}\{e^{itS_{r,n}}\}^2\right|$$
$$\leq (2+|t|)\mathcal{E}^{\xi}_{h}(n) + 16\mathcal{E}^{I}_{k,h}(n) + 14\mathcal{E}^{S}_{k,h}(n) + 20.4^{-k} + 20|t|4^{-k} \quad \text{(V-7)} .$$

---

**Proof :** We know that $|e^{ib}-e^{ia}|$ for all (a,b) $\mathbb{R}^2$ . Then,

$$\left|\mathbb{E}\{e^{it\psi_{r-1}(N).S_{r-1,n}/\psi_r(N)}\} - \mathbb{E}\{e^{itS_{r,n}+itS'_{r,n}}\}\right|$$
$$\leq \mathbb{E}\left\{|e^{itS_{r,n}+itS'_{r,n}}|\ |e^{it\xi_{r,n}} - 1|\right\} \leq \mathbb{E}\left\{|e^{it\xi_{r,n}} - e^{0}|\right\}$$
$$\leq \mathbb{E}_{|\xi_{r,n}|>\mathcal{E}^{\xi}_{h}(n)}\{2\} + \mathbb{E}_{|\xi_{r,n}|\leq\mathcal{E}^{\xi}_{h}(n)}\left\{|t.\xi_{r,n} - 0|\right\} \leq (2+|t|)\mathcal{E}^{\xi}_{h}(n) \quad .$$

Then, it is enough to use this inequality and (V-5) and (V-6) about

$$\left|\mathbb{E}\{e^{it\psi_{r-1}(N).S_{r-1,n}/\psi_r(N)}\} - \mathbb{E}\{e^{itS_{r,n}+itS'_{r,n}}\}\right|$$
$$\left|\mathbb{E}\{e^{itS_{r,n}+itS'_{r,n}}\} - \mathbb{E}\{e^{itS_{r,n}}\}\ \mathbb{E}\{^{itS'_{r,n}}\}\right|$$

and

$$\left|\mathbb{E}\{e^{itS_{r,n}}\}\left(\mathbb{E}\{e^{itS'_{r,n}}\} - \mathbb{E}\{e^{itS_{r,n}}\}\right)\right| \quad .$$

---

**Lemma 5-8 :** For all t $\mathbb{R}$ all k and h $\mathbb{N}$ , and all r=1,2,...,h,

$$\left|\mathbb{E}\{e^{itS_{r-1,n}}\} - \left(\mathbb{E}\{e^{itS_{r,n}/2^t}\}\right)^2\right|$$
$$\leq b_k|t|\mathcal{E}^{\psi}_{\tau,h}(n) + (2+2^{-\tau}|t|)\mathcal{E}^{\xi}_{h}(n) + 16\mathcal{E}^{I}_{k,h}(n) + 14\mathcal{E}^{S}_{k,h}(n) + 4^{-k}(22 + 20|t|2^{-\tau}) \quad \text{(V-8)} .$$

---

**Proof :** By (V-1), one can write

$$\left|\mathbb{E}\{e^{itS_{r-1,n}}\} - \mathbb{E}\{e^{it2^{-\tau}\psi_{r-1}(N).S_{r-1,n}/\psi_r(N)}\} - \right|$$

$$\leq \mathbb{E}_{|S_{r-1,n}|\leq b_k}\left\{ \left| e^{itS_{r-1,n}\}} - e^{it2^{-\tau}\psi_{r-1}(N).S_{r-1,n}/\psi_r(N)} \right| \right\} + 2.4^{-k}$$

$$\leq \mathbb{E}_{|S_{r-1,n}|\leq b_k}\left\{ \left| tS_{r-1,n} - t.2^{-\tau}\psi_{r-1}(N).S_{r-1,n}/\psi_r(N) \right| \right\} + 2.4^{-k}$$

$$\leq t.b_k \mathcal{E}^{\psi}_{\tau,h}(n) \quad + 2.4^{-k} .$$

The result follows from (V-7).

---

**Lemma 5-9 :** For all $r \ \mathbb{R}$ , we set $\delta_r = 2^{-r\tau}$ . For all $r \ \mathbb{R}$ and all k and h $\mathbb{N}$ ,

$$\left| \mathbb{E}\{\exp(itS_{0,n})\} - \mathbb{E}\{\exp(\delta_h itS_{h,n})\}^{2^h} \right|$$

$$\leq 2^h\Big(b_{k+h}|t|\mathcal{E}^{\psi}_{\tau,h}(n) + (2+|t|)\mathcal{E}^{\xi}_{h}(n) + 16\mathcal{E}^{I}_{k+h,h}(n) + 14\mathcal{E}^{S}_{k+h,h}(n)\Big)+4^{-k}(44 + 40|t|) \quad \text{(V-9)} .$$

---

**Proof :** Let k’≤k. Clearly, $\mathfrak{T}_k \subset \mathfrak{T}_{k'}$ . Then, $\mathcal{E}^{S}_{k,h}(n) \leq \mathcal{E}^{S}_{k',h}(n)$ and $\mathcal{E}^{I}_{k,h}(n) \leq \mathcal{E}^{I}_{k',h}(n)$.
Therefore, by using (V-8),

$$\left| \mathbb{E}\{\exp(itS_{0,n})\} - \mathbb{E}\{\exp(\delta_h itS_{h,n})\}^{2^h} \right|$$

$$\leq \left| \sum_{r=0}^{h-1} \Big[\mathbb{E}\{\exp(\delta_r itS_{r,n})\}^{2^r} - \mathbb{E}\{\exp(\delta_{r+1} itS_{r+1,n})\}^{2^{r+1}}\Big] \right|$$

$$\leq \sum_{r=0}^{h-1} 2^r \left| \mathbb{E}\{\exp(\delta_r itS_{r,n})\} - \mathbb{E}\{\exp(\delta_{r+1} itS_{r+1,n})\}^{2^r} \right|$$

$$\leq \sum_{r=0}^{h-1} 2^r\, b_{k+r}2^{-r\tau}|t|\mathcal{E}^{\psi}_{\tau,h}(n) + \sum_{r=0}^{h-1} 2^r(2+2^{-(r+1)\tau}|t|)\mathcal{E}^{\xi}_{h}(n)$$
$$+ 16\sum_{r=0}^{h-1} 2^r\mathcal{E}^{I}_{k+r,h}(n) + 14\sum_{r=0}^{h-1} 2^r\mathcal{E}^{S}_{k+r,h}(n) + 22\sum_{r=0}^{h-1} 2^{-2k-r} + 20\sum_{r=0}^{h-1} 2^{-2k-r}2^{-(r+1)\tau}|t|$$

$$\leq |t|\sum_{r=0}^{h-1} 2^r b_{k+h}\mathcal{E}^{\psi}_{\tau,h}(n) + (2+|t|)\sum_{r=0}^{h-1} 2^r\mathcal{E}^{\xi}_{h}(n)$$
$$+ 16\sum_{r=0}^{h-1} 2^r\mathcal{E}^{I}_{k+h,h}(n) + 14\sum_{r=0}^{h-1} 2^r\mathcal{E}^{S}_{k+h,h}(n) + 4^{-k}(44+40|t|) .$$

---

**Lemma 5-10 :** Let $h_n$ $\mathbb{N}^*$ be a non-decreasing sequence. Then, for all $t$ $\mathbb{R}$,

$$\left| \mathbb{E}\{\exp(itS_{0,n})\} - \mathbb{E}\{\exp(\delta_{h_n} itS_{h_n,n})\}^{2^{h_n}} \right|$$

$$\leq 2^{h_n}\Big(b_{2h_n}|t|\mathcal{E}^{\psi}_{\tau,h_n}(n) + (2+|t|)\mathcal{E}^{\xi}_{h_n}(n) + 16\mathcal{E}^{I}_{2h_n,h_n}(n) + 14\mathcal{E}^{S}_{2h_n,h_n}(n)\Big)$$

$$+4^{-h_n}(44 + 40|t|) \quad \text{(V-10)} .$$

---

## V-3 : Some propositions

By using the previous lemmas, we obtain the following results.

**Proposition 5-11:** Assume that, for all $k$ $\mathbb{N}^*$ , $\mathcal{E}^{\tau}_{k}(n) \to 0$ as $n \to \infty$ and that, for all $k$ $\mathbb{N}^*$, all $\tau>0$, all $n$ $\mathbb{N}^*$ , $\mathbf{P}\{S_{r,n} \ \mathfrak{K}_k\} \leq 4^{-k}$ .

Then, $S_{0,n}$ converges in distribution to a random variable L if and only if there exists a non decreasing sequence $k_n$ $\mathbb{R}^*$ , $k_n \to \infty$ as $n \to \infty$ , such that the two following assertions hold.

$$2^{k_n}\mathcal{E}^{\tau}_{k_n}(n) \to 0 \text{ as } n \to \infty \quad \text{(V-11)} ,$$

$$2^{-\tau k_n}\Big(\sum\nolimits^{n}_{1} + \sum\nolimits^{n}_{2} + \ldots\ldots. + \sum\nolimits^{n}_{2^{k_n}}\Big) \xrightarrow{d} L \quad \text{(V-12)} ,$$

where the $\sum^{n}_{j}$'s , $j=1,2,....,2^{k_n}$, are independent random variables which have the same distribution as $S_{k_n,n}$ .

---

**Proof :** Assume that (V-11) and (V-12) hold.

By (V-12), $\mathbb{E}\{\exp(itL)\} - \mathbb{E}\{\exp(\delta_{k_n} itS_{k_n,n})\}^{2^{k_n}} \to 0$ as $n \to \infty$.

By (V-10) and (V-11), $\mathbb{E}\{\exp(itS_{0,n})\} - \mathbb{E}\{\exp(\delta_{k_n} itS_{k_n,n})\}^{2^{k_n}} \to 0$ as $n \to \infty$.

We deduce that $\mathbb{E}\{\exp(itS_{0,n})\} \to \mathbb{E}\{\exp(itL)\}$ as $n \to \infty$.

Assume that $S_{0,n}$ converges to L. It is always possible to choose an increasing sequence $\{k_n\}$ which increases enough slowly in order that (V-11) holds. For a such sequence, by (V-10), $\mathbb{E}\{\exp(\delta_{k_n} itS_{k_n,n})\}^{2^{k_n}} \to \mathbb{E}\{\exp(itL)\}$ as $n \to \infty$.

---

Now, by (V-9) the following result holds.

**Corollary 5-12 :** Assume that all the assumptions of proposition 5-11 hold. Then, for all r $\mathbb{N}$, $\mathbb{E}\{\exp(2^{-r\tau} itS_{r,n})\}^{2^r} \to \mathbb{E}\{\exp(itL)\}$ as $n \to \infty$.

---

In particular, if $\mathbb{E}\{\exp(itL)\} = \exp(-t^2/2)$ and $\tau=1/2$, $S_{r,n} \xrightarrow{d} N(0,1)$ for all r $\mathbb{N}$. Then, we can apply proposition 5-11 for the convergence to the normal distribution. In this cas, we choose the following assumptions.

**Assumptions 5-13 :** We assume that, for all n $\mathbb{N}$, $\mathbb{E}\{X_{n,t}\} = 0$, $\mathbb{E}\{(X_{n,t})^2\} <\infty$ and $\sigma_r(n)^2 = \Psi_r(n)^2 = \mathbb{E}\{(X_{n,1}+X_{n,2}+.....+X_{n,u^r(n)})^2\} <\infty$. We assume that, for all k $\mathbb{N}$, $\mathfrak{Z}_k$ is the samallest subset of $\mathbb{N}$ such that $[-2^k,2^k]\subset\mathfrak{K}_k$. We assume also that, for all k and $\mathbb{N}$, $\varepsilon^I_{k,h}(n)\to 0$, $\varepsilon^S_{k,h}(n)\to 0$ and $\varepsilon^\xi_h(n)\to 0$ as $n \to \infty$.

---

Then, one generalizes the theorem 1-3 of [13]. In this case, $\tau=1/2$.

**Corollary 5-14 :** Assume that the assumptions of 5-13 hold. Then, $S_{r,n} \xrightarrow{d} N(0,1)$ for all r $\mathbb{N}$, if and only if there exists a non-decreasing sequence $k_n$ $\mathbb{N}^*$, $k_n\to\infty$ as $n \to \infty$, such that (V-11) and (V-13) hold for $\tau=1/2$, where (V-13) is the Lindeberg condition :

$$\forall d>0,\ \mathbb{E}_{|S_{k_n,n}|>d2^{k_n/2}}\{(S_{k_n,n})^2\} \to 0 \text{ as } n \to\infty \quad \text{(V-13)}.$$

---

**Proof :** By the Bienaymé-Tschébyscheff Inequality, $\mathbf{P}\{S_{r,n}\ \mathfrak{K}_k)\} \leq 4^{-k}$. Moreover, by the theorem of page 103 of [39], (V-12) is equivalent to (V-13) if $L\sim N(0,1)$ and $\tau=1/2$ .

Assume that (V-11) and (V-13) hold. By proposition 5-11, $S_{0,n} \xrightarrow{d} N(0,1)$. By (V-9), $S_{r,n} \xrightarrow{d} N(0,1)$ for all r $\mathbb{N}$ .

Assume that $S_{r,n} \xrightarrow{d} N(0,1)$ for all r $\mathbb{N}$ . By (V-7), $\mathbb{E}\{e^{it\sigma_{r-1}(N).S_{r-1,n}/\sigma_r(N)}\} \to \exp(-t^2/2)$ as $n \to\infty$ . Therefore, $\sigma_{r-1}(N)S_{r-1,n}/\sigma_r(N) \xrightarrow{d} L'\sim N(0,2)$ . Therefore,

$$\text{for all r } \mathbb{N}^*,\quad \sigma_{r-1}(N)/\sigma_r(N) \to \sqrt{2} \text{ as } n \to\infty \text{ (V-14)} .$$

Then, we can use theorem 5-11 with τ=1/2 .

---

Now, one can always build a proper sequence $k_n$ if the Lindeberg Condition holds.

**Corollary 5-15 :** Assume that the assumptions of 5-13 hold. Then, $S_{r,n} \xrightarrow{d} N(0,1)$ for all $r \in \mathbb{N}$, if and only if (V-14) holds and there exists an increasing sequence $N_k$ such that for all $n \geq N_k$,

$$\forall d>0, \quad \mathbb{E}_{|S_{r,n}|>2^{k/4}}\left\{(S_{r,n})^2\right\} \leq \Delta'_k \text{ for all } r \leq k \qquad \text{(V-15)}$$

where $\Delta'_k \to 0$ as $k \to \infty$ .

---

**Proof :** Assume that $N_k$ is strictly increasing. We define $h_n$ by $h_n=1$ if $n<N_2$ and $h_n=s$ if $N_s \leq n < N_{s+1}$.

Then, for all $m \geq N_{h_n}$ , $\mathbb{E}_{|S_{r,m}|>2^{h_n/4}}\left\{(S_{r,m})^2\right\} \leq \Delta'_{h_n}$ for all $r \leq h_n$ .

Let $k_n$ be a sequence such that $k_n \leq h_n$ and $k_n \to \infty$ . Because $h_{n+1} = h_n$ or $h_{n+1} = h_n+1$ , there exists t such that $k_n = h_t$ . Therefore, for all $m \geq N_{h_t}$ , $\mathbb{E}_{|S_{r,m}|>2^{k_n/4}}\left\{(S_{r,m})^2\right\} = \mathbb{E}_{|S_{r,m}|>2^{h_t/4}}\left\{(S_{r,m})^2\right\} \leq \Delta'_{h_t}$ for any $r \leq h_t = k_n$ . Moreover, $n \geq N_{hn_t} \geq N_{h_t}$ . Then, one can apply the previous inequality with m=n and $r \leq k_n$ : $\mathbb{E}_{|S_{k_n,n}|>2^{k_n/4}}\left\{(S_{k_n,n})^2\right\} \leq \Delta'_{k_n}$ .

Therefore, (V-13) holds for all sequences $k_n \leq h_n$ .

Moreover, one can choose $k_n$ which increase slowly in order that (V-11) holds.
On the other hand, the proof of the necessary condition is classical (e.g. cf [13], p 339).

---

### V-4 : Proof of theorem 3-1.

In order to prove theorem 3-1, the following lemma is needed.

**Lemma 5-16 :** Assume that the assumptions of 5-13 hold. Moreover, we assume that, for all $r \ \mathbb{N}$, $\mathbb{E}\{(S'_{r,n})^2\} - \mathbb{E}\{(S_{r,n})^2\} = \varepsilon_2^r(n) \to 0$ and $\mathbb{E}\{(\xi_{r,n})^2\} \to 0$ as $n \to \infty$ , and that

$$\mathbb{E}_{|S_{r,n}|>k}\{(S_{r,n})^2\} \leq \Delta_k^r \to 0 \text{ as } k \to \infty \quad \text{(V-16)} .$$

Then, (V-14) holds.

---

**Proof :** One can assume that $\pi_k^r = \Delta_{2k}^r$ is non-increasing. We know that $b_k \leq 2^{k+1}$ and $\mathbf{P}\{S_{r,n} \ \mathcal{K}_k)\} \leq 4^{-k}$.

There exists $I_k^2 = \sum_{i \ \mathfrak{J}_k} \beta_i^k ¶_{\mathbf{J}_{k,i}}$ , $|\beta_i^k| \leq 4^{k+1}$ such that $|s^2 - I_k^2(s)| \leq (2.2^{k+1})4^{-k}$ if $s \ \mathcal{K}_k$.

Then, by (V-3), for r>0,

$$\left|\mathbb{E}_{S_{r,n} \ \mathcal{K}_k}\{(S_{r,n})^2\} - \mathbb{E}_{S'_{r,n} \ \mathcal{K}_k}\{(S'_{r,n})^2\}\right|$$

$$\leq \left| 1 - \mathbb{E}_{S_{r,n} \ \mathcal{K}_k}\{(S_{r,n})^2\} - \left(1 + \varepsilon_2^r(n) - \mathbb{E}_{S'_{r,n} \ \mathcal{K}_k}\{(S'_{r,n})^2\}\right)\right|$$

$$\leq \left|\mathbb{E}_{S_{r,n} \ \mathcal{K}_k}\{I_k^2(S_{r,n})\} - \mathbb{E}_{S'_{r,n} \ \mathcal{K}_k}\{I_k^2(S'_{r,n})^2\}\right| + 8.2^{-k} + |\varepsilon_2^r(n)|$$

$$\leq 2.4^{k+1}\varepsilon_{k,r}^S(n) + 8.2^{-k} + |\varepsilon_2^r(n)| \ .$$

Then, for all $k \ \mathbb{N}$, $\left|\mathbb{E}_{S'_{r,n} \ \mathcal{K}_k}\{(S'_{r,n})^2\}\right| \leq \pi_k^r + 2.4^{k+1}\varepsilon_{k,r}^S(n) + 8.2^{-k} + |\varepsilon_2^r n)|$ .

Let $N_k^r$ such that, for all $n \geq N_k^r$, $\varepsilon_2^r(n) \leq 2^{-k}$ and $\varepsilon_{k,r}^S(n) \leq 8^{-k}$ . We assume $N_k^r$ increasing. We define $\mu_k^r > 0$ by $(\mu_k^r)^2 = \max\left\{\pi_k^r + 17.2^{-k} \ , \ \operatorname{Max}_{n<N_k^r} \left(\mathbb{E}_{S'_{r,n} \ \mathcal{K}_k}\{(S'_{r,n})^2\}\right) \right\}$ .

Then, the following lemma is needed.

**Lemma 5-17 :** Under the previous assumptions

$$\mu_k^r \to 0 \text{ as } k \to \infty \text{ (V-17) .}$$

---

**Proof :** let $\gamma>0$. There exists K such that $\pi_K^r + 17.2^{-K} \leq \gamma$. Then, there exists K'>0 such that

$$\operatorname*{Max}_{n<N_K^r} \left( \mathbb{E}_{S'_{r,n}\ \mathfrak{K}_{K'}} \left\{ (S'_{r,n})^2 \right\} \right) \leq \gamma.$$

Let $k \geq \max(K,K')$ . Then,

$$\operatorname*{Max}_{n<N_k^r} \left( \mathbb{E}_{S'_{r,n}\ \mathfrak{K}_{k}} \left\{ (S'_{r,n})^2 \right\} \right)$$

$$\leq \max\left\{ \operatorname*{Max}_{n<N_K^r} \left( \mathbb{E}_{S'_{r,n}\ \mathfrak{K}_{K'}} \left\{ (S'_{r,n})^2 \right\} \right) ,\quad \operatorname*{Max}_{N_K^r \leq n < N_k^r} \left( \mathbb{E}_{S'_{r,n}\ \mathfrak{K}_{K}} \left\{ (S'_{r,n})^2 \right\} \right) \right\}$$

$$\leq \gamma .$$

We deduce (V-17).

---

Clearly, $\mathbb{E}_{S'_{r,n}\ \mathfrak{K}_{k}} \left\{ (S'_{r,n})^2 \right\} \leq (\mu_k^r)^2$ and $\mathbb{E}_{S_{r,n}\ \mathfrak{K}_{k}} \left\{ (S_{r,n})^2 \right\} \leq (\mu_k^r)^2$ .

Moreover, there exists $B_r>0$ such that $\mathbb{E}\left\{ |S'_{r,n}| \right\} \leq B_r$ and $\mathbb{E}\left\{ |S'_{r,n}|^2 \right\} \leq (B_r)^2$ .

Then, by Schwartz Inequality, $\left| \mathbb{E}_{S_{r,n}\ \mathfrak{K}_{k}} \left\{ S_{r,n} S'_{r,n} \right\} \right| \leq B_r\, \mu_k^r$ , for example.

Then,

$$\left|\mathbb{E}\left\{S_{r,n}S'_{r,n}\right\} - \mathbb{E}\left\{S_{r,n}\right\}\mathbb{E}\left\{S'_{r,n}\right\}\right|$$

$$\leq \left|\mathbb{E}_{S_{r,n}\ \mathfrak{K}_k}\left\{S_{r,n}S'_{r,n}\right\} - \mathbb{E}_{S_{r,n}\ \mathfrak{K}_k}\left\{S_{r,n}\right\}\mathbb{E}\left\{S'_{r,n}\right\}\right|$$
$$+ \left|\mathbb{E}_{S_{r,n}\ \mathfrak{K}_k}\left\{S_{r,n}S'_{r,n}\right\}\right| + B_r\left|\mathbb{E}_{S_{r,n}\ \mathfrak{K}_k}\left\{S_{r,n}\right\}\right|$$

$$\leq \left|\mathbb{E}_{S_{r,n}\ \mathfrak{K}_k}\left\{S_{r,n}S'_{r,n}\right\} - \mathbb{E}_{S_{r,n}\ \mathfrak{K}_k}\left\{S_{r,n}\right\}\mathbb{E}\left\{S'_{r,n}\right\}\right| + 2\,B_r\,\mu_k^r$$

$$\leq \left|\mathbb{E}_{\{S_{r,n}\ \mathfrak{K}_k\}\cap\{S'_{r,n}\ \mathfrak{K}_k\}}\left\{S_{r,n}S'_{r,n}\right\} - \mathbb{E}_{S_{r,n}\ \mathfrak{K}_k}\left\{S_{r,n}\right\}\mathbb{E}_{S'_{r,n}\ \mathfrak{K}_k}\left\{S'_{r,n}\right\}\right|$$
$$+ 2B_r\mu_k^r + 2\mu_k^r\ .$$

Moreover, there exists $I_k = \sum_{i\ \mathfrak{J}_k} \beta_i^k \P_{k,i}$ , $|\alpha_i^k| \leq 2^{k+1}$ , such that $|s - I_k(s)| \leq 4^{-k}$ and $|I_k(s)| \leq s$ if $s\ \ \mathfrak{K}_k$.

Then,

$$\left|\mathbb{E}_{\{S_{r,n}\ \mathfrak{K}_k\}\cap\{S'_{r,n}\ \mathfrak{K}_k\}}\left\{S_{r,n}S'_{r,n}\right\} - \mathbb{E}_{S_{r,n}\ \mathfrak{K}_k}\left\{S_{r,n}\right\}\mathbb{E}_{S'_{r,n}\ \mathfrak{K}_k}\left\{S'_{r,n}\right\}\right|$$

$$\leq \left|\mathbb{E}_{\{S_{r,n}\ \mathfrak{K}_k\}\cap\{S'_{r,n}\ \mathfrak{K}_k\}}\left\{I_k(S_{r,n})S'_{r,n}\right\} - \mathbb{E}_{S_{r,n}\ \mathfrak{K}_k}\left\{I_k(S_{r,n})\right\}\mathbb{E}_{S'_{r,n}\ \mathfrak{K}_k}\left\{S'_{r,n}\right\}\right|$$
$$+ 4^{-k}\left|\mathbb{E}_{\{S_{r,n}\ \mathfrak{K}_k\}\cap\{S'_{r,n}\ \mathfrak{K}_k\}}\left\{|S'_{r,n}|\right\}\right| + 4^{-k}\left|\mathbb{E}_{S'_{r,n}\ \mathfrak{K}_k}\left\{S'_{r,n}\right\}\right|$$

$$\leq \left|\mathbb{E}_{\{S_{r,n}\ \mathfrak{K}_k\}\cap\{S'_{r,n}\ \mathfrak{K}_k\}}\left\{I_k(S_{r,n})S'_{r,n}\right\} - \mathbb{E}_{S_{r,n}\ \mathfrak{K}_k}\left\{I_k(S_{r,n})\right\}\mathbb{E}_{S'_{r,n}\ \mathfrak{K}_k}\left\{S'_{r,n}\right\}\right|$$
$$+ 2B_r\,4^{-k}$$

$$\leq \Big| \mathbb{E}_{\{S_{r,n} \ \mathfrak{K}_k\} \cap \{S'_{r,n} \ \mathfrak{K}_k\}} \big\{ I_k(S_{r,n}) I_k(S'_{r,n}) \big\}$$

$$- \mathbb{E}_{S_{r,n} \ \mathfrak{K}_k} \big\{ I_k(S_{r,n}) \big\} \mathbb{E}_{S'_{r,n} \ \mathfrak{K}_k} \big\{ I_k(S'_{r,n}) \big\} \Big|$$

$$+ 4^{-k} \mathbb{E}_{\{S_{r,n} \ \mathfrak{K}_k\} \cap \{S'_{r,n} \ \mathfrak{K}_k\}} \big\{ |I_k(S_{r,n})| \big\} + 4^{-k} \mathbb{E}_{S_{r,n} \ \mathfrak{K}_k} \big\{ |I_k(S_{r,n})| \big\} + 2B_r 4^{-k}$$

$$\leq \ 4.4^{k+1} \varepsilon^{I}_{k,r}(n) \ + \ 2(1+B_r) 4^{-k} \ .$$

These inequalities hold for all k. Therefore, by (V-17), $\mathbb{E}\{S_{r,n} S'_{r,n}\}$ converges to 0. Moreover, $\mathbb{E}\{S_{r,n} \xi_{r,n}\}^2 \ \leq \ \mathbb{E}\{(S_{r,n})^2\} \mathbb{E}\{(\xi_{r,n})^2\}$ which converges to 0. Therefore, $\mathbb{E}\{(S_{r,n} + \xi_{r,n} + S'_{r,n})^2\}$ converges to 2. Then, (V-14) holds.

---

**Lemma 5-18 :** Under the assumptions of definition 3-1, $n/(u(n) \to 2$ as $n \to \infty$. Moreover, if $\mathbb{E}\{(S'_u)^2\} - \mathbb{E}\{(S_u)^2\} \to 0$ and $\mathbb{E}\{(\xi_u)^2\} \to 0$ as $n \to \infty$, $\sigma(n)^2/\sigma(u)^2 \to 2$ and $\sigma(n) \to \infty$ as $n \to \infty$ .

---

**Proof :** by notations 2-1, if $n>1$, $3(u+1) \geq 2(u+1) + \kappa(u+1) > n > 2u + \kappa(u)$. Therefore, $u(n) \to \infty$ as $n \to \infty$ and $2 + 2/u + [\kappa(u+1)/(u+1)]\ [(u+1)/u] \ > n/u \geq 2 + \kappa(u)/u$. Then, by notations 2-1, $n/(u(n) \to 2$ as $n \to \infty$.

Assume that $\mathbb{E}\{(S'_u)^2\} - \mathbb{E}\{(S_u)^2\} \to 0$ and $\mathbb{E}\{(\xi_u)^2\} \to 0$ . Then, $1 = \mathbb{E}\{(S_n)^2\} = [\sigma(u)^2/\sigma(n)^2] \mathbb{E}\{(S_u + \xi_u + S'_u)^2\} = [\sigma(u)^2/\sigma(n)^2][\mathbb{E}\{(S_u)^2\} + o(1)] = [\sigma(u)^2/\sigma(n)^2][2 + o(1)]$ . Then, $\sigma(n)^2/\sigma(u)^2 \to 2$ as $n \to \infty$ .

Therefore, there exists $N_0$ such that, for $n \geq N_0$, $\sigma(n)^2/\sigma(u)^2 \geq 5/3$ and $7/3 \geq n/u(n)$ . Then, if $n \geq (7/3)^k N_0$, $u^k(n) \geq N_0$ . Therefore, there exists $k' \geq k$ such that $u^{k'-1}(n) \geq N_0$ and $u^{k'}(n) < N_0$ . Let $M = \min\{\sigma(n)^2 |\ n < N_0\}$ . Then, $\sigma(n)^2 \geq (5/3)^k M$. Because, $\sigma(n) > 0$, $\sigma(n) \to \infty$ as $n \to \infty$ .

---

**5-19 : Proof of theorem 3-1.** Let $M_r$ be an increasing sequence properly choosen : e.g. one can assume $u^{r_n} \to \infty$ where $r_n = r$ if $M_r < n \leq M_{r+1}$ . We set $V(n) = \sigma(n) \upsilon_n$ , $V'(n) = \sigma(n) \upsilon'_n$ and $N(n) = n$.

If $n \leq M_1$, we define $X_{n,t}$ by $X_{n,t} = X_t$ if $t<n$ and $X_{n,n} = X_n + V_n$ .

If $M_1<n\leq M_2$, we set $X_{n,t} = X_t$ if $t\neq u,u+1,n$ and $X_{n,u} = X_u + V_u$ , $X_{n,u+1} = X_{u+1}+V_n-V_u-V'_u$ , $X_{n,n} = X_n+V'_u$ .

If $M_h<n\leq M_{h+1}$, we set $X_{n,t} = X_t$ if $t\neq u^r(n),u^r(n)+1$, for $r=0,1,...,h$ and $X_{n,u^h(n)} = X_{u^h(n)}+V_{u^h(n)}$ , $X_{n,u^r(n)} = X_{u^r(n)}+V'_{u^{r+1}(n)}$ for $r=0,1,2,.....,h-1$ , $X_{n,u^r(n)+1} = X_{u^r(n)+1} + V_{u^{r-1}(n)} - V_{u^r(n)} - V'_{u^r(n)}$ for $r=1,2,...,h$ (cf also [37]).

Then, for $n>M_r$ ,

$$\sigma_r(n)S_{r,n} = X_1+X_2+.....+ X_{u^r(n)}+V_{u^r(n)} = \sigma_r(u^r(n))( S_{u^r(n)}+\upsilon_{u^r(n)}) ,$$

$$\sigma_r(n)S'_{r,n} = X_{u^r(n)+t^r(n)+1} + X_{u^r(n)+t^r(n)+2} +..........+ X_{u^{r-1}(n)} + V'_{u^r(n)}$$
$$= \sigma_r(u^r(n))(S'_{u^r(n)}+\upsilon'_{u^r(n)}) ,$$

$$\sigma_r(n)\xi_{r,n} = X_{u^r(n)+1} + X_{u^r(n)+2} +..........+ X_{u^r(n)+t^r(n)} + V_{u^{r-1}(n)} - V_{u^r(n)} + V'_{u^r(n)} .$$

Therefore, $\sigma_r(N)^2 = \mathbb{E}\{(\sigma_r(N)S_{r,n})^2\} = \sigma(u^r(n))^2\mathbb{E}\{(S_{u^r(n)}+\upsilon_{u^r(n)})^2\} = \sigma(u^r(n))^2\sigma^*(u^r(n))^2$ .

Therefore, $S_{r,n} = \sigma^*(u^r(n))^{-1}(S_{u^r(n)}+\upsilon_{u^r(n)})$ and $S'_{r,n}= \sigma^*(u^r(n))^{-1}(S'_{u^r(n)}+\upsilon'_{u^r(n)})$.

Then, $\varepsilon^I_{k,h}(n) \to 0$ and $\varepsilon^S_{k,h}(n) \to 0$ .

By our assumptions, $\sigma^*(u^r(n))^2 \to 1$ .

Then, $\mathbb{E}\{(S_{r,n})^2\}-\mathbb{E}\{(S'_{r,n})^2\} \to 0$ .

Moreover, by lemma 5-18, $\sigma(u^r(n))^2/\sigma(u^{r-1}(n))^2 \to 2$ . Then, $\mathbb{E}\{(\xi_{r,n})^2\} \to 0$ Therefore, one can assume $\varepsilon^\xi_h(n) \to 0$ .

On the other hand, by our assumptions, there exists $\Delta_k^r \to 0$ as $k \to \infty$ such that

$$\mathbb{E}_{|S_{r,n}|\geq k}\left\{(S_{r,n})^2\right\} = \mathbb{E}_{\left|\sigma^*(u^r(n))^{-1}(S_{u^r(n)}+\upsilon_{u^r(n)})\right|\geq k}\left\{\sigma^*(u^r(n))^{-2}\left(S_{u^r(n)}+\upsilon_{u^r(n)}\right)^2\right\} \leq \Delta_k^r \ .$$

Then, (V-14) holds.

Then all the assumptions of corollary 5-15 are satisfied by $\{X_n\}$. Therefore, $S_n \xrightarrow{d} N(0,1)$ .

---